\documentclass[12pt]{amsart}

\usepackage[latin2]{inputenc}

\usepackage{amsmath, amssymb}
\newtheorem{theorem}{Theorem}
\newtheorem{conjecture}[theorem]{Conjecture}
\newtheorem{corollary}[theorem]{Corollary}
\newtheorem{proposition}[theorem]{Proposition}
\newtheorem{lemma}[theorem]{Lemma}

\newtheorem{definition}{Definition}

\def\B{\mathcal{B}}

\def\G{\mathcal{G}}

\def\bth{\begin{theorem}}
\def\eth{\end{theorem}}

\def\bc{\begin{corollary}}
\def\ec{\end{corollary}}

\def\bcj{\begin{conjecture}}
\def\ecj{\end{conjecture}}

\newcommand{\bizveg}{{\hfill$\Box$}}

\title[On the Complexity of Chooser-Picker Positional Games]{On the Complexity of Chooser-Picker Positional Games}

\author[Csernenszky]{Andr\'as Csernenszky}
\address{University of Szeged, Department of Computer Science}
\email{csernenszkya@gmail.com}

\author[Martin]{Ryan R. Martin}
\thanks{This author's research partially supported by NSF grant DMS-0901008.}
\address{Iowa State University, Department of Mathematics}
\email{rymartin@iastate.edu}

\author[Pluh\'ar]{Andr\'as Pluh\'ar}
\thanks{The research was partially supported by OTKA grants T049398 and K76099.}
\address{University of Szeged, Department of Computer Science}
\email{pluhar@inf.u-szeged.hu}


\begin{document}

\begin{abstract}
Two new versions of the so-called Maker-Breaker Positional Games
are defined by J\'ozsef Beck in [{\em Combinatorica} {\bf 22}(2) (2002) 169--216]. He defines two players, Picker and Chooser.  In each round, Picker takes a pair of elements not already selected and Chooser keeps one and returns the other to Picker. In the Picker-Chooser version
Picker plays as Maker and Chooser plays as Breaker, while the roles are swapped in the Chooser-Picker version. The outcome of these games is sometimes very similar to that of the traditional Maker-Breaker games. Here we show that both Picker-Chooser and Chooser-Picker games are NP-hard, which gives support to the paradigm that the games behave similarly while being quite different in definition. We also investigate the pairing strategies for Maker-Breaker games, and apply these results to the game called ``Snaky.''
\end{abstract}

\maketitle

\section{Introduction}
Let us start by defining Positional Games in general. Given an arbitrary hypergraph ${\mathcal H}$ with vertex set $V(\mathcal H)$ and edge set $E(\mathcal H)$, we write ${\mathcal H}=(V, E)$ and the first and second players take elements of $V$ in turns. The goal of a player designated as a \textit{Maker} is to take every element of some edge $A\in E$.  In the Maker-Maker version of the game, the player who is first to take all elements of some edge $A\in E$ wins the game. There are a number of beautiful and difficult theorems about Maker-Maker games, for more details Berlekamp, Conway and Guy:~\cite{B-C-G} or Beck:~\cite{beckcpc05} or~\cite{BB}.

The so-called Maker-Breaker version of a Positional Game on a hypergraph ${\mathcal H}=(V, E)$ was also investigated from the very beginning. Here each player takes an unselected element in turn.  One player is still Maker and wins by taking every element of some $A\in E$.  The other player is \textit{Breaker} and wins if he can take at least one vertex of every edge in $E$.  Clearly, exactly one of Maker and Breaker can win this game. In several cases Maker-Breaker games are more tractable than Maker-Maker games. On the other hand, these versions are closely related, since if Breaker wins as a second player then the Maker-Maker game is a draw; that is, the second player can ensure that the first does not have a winning strategy.  On the other hand, if the first player has a winning strategy for the Maker-Maker game, then Maker also wins the Maker-Breaker version. This connection gives rise to very useful applications, see~\cite{B2, Beck3, Beck2, Zetters, H-J}.

In order to understand the so-called clique games, which are very difficult, Beck introduced the Picker-Chooser and the Chooser-Picker version of Maker-Breaker games in~\cite{Beck5}.


\begin{definition}
The positional game players \textbf{Picker} and \textbf{Chooser} are as follows: Picker takes a pair of elements, neither of which had been selected previously, and Chooser keeps one of these elements and gives the other to Picker. The designation ``\textbf{Picker-Chooser}'' indicates that Picker plays as Maker (that is, wins by taking all the elements in some edge) and Chooser plays as Breaker (that is, wins by taking at least one element from each edge). The roles are swapped in the ``\textbf{Chooser-Picker}'' version, in which Chooser plays as Maker and Picker plays as Breaker. If $|V|$ is odd, then the last element goes to Chooser.
\end{definition}

Beck demonstrated in several cases that Picker may easily win the Picker-Chooser game if Maker wins the corresponding Maker-Breaker game~\cite{Beck5,BB}.

A similar phenomenon must also hold for certain Chooser-Picker games; that is, Picker is better off in the Chooser-Picker version than Breaker is in the corresponding Maker-Breaker game. In fact, there is a duality between Maker-Breaker games in which Maker wins and ones in which Breaker wins. Let ${\mathcal H^*}=(V, E^*)$ be the \textit{transversal hypergraph} of ${\mathcal H}=(V, E)$. That is, $E^*$ consists of those minimal sets $B\subset V$ such that for all $A\in E$, $A\cap B\neq\emptyset$. Note that Breaker as a first [second] player wins the Maker-Breaker game on $(V, E)$ iff Maker as a first [second] player wins the Maker-Breaker $(V, E^*)$.

The general form of Beck's conjecture is spelled out by Csernenszky, et al.~\cite{CMP}.

\bcj \cite{CMP} \label{Bc}
Picker wins a Picker-Chooser [Chooser-Picker] game on $(V, E)$ if Maker
[Breaker] as a second player wins the corresponding Maker-Breaker game.
\ecj

Note that this generalized Beck's conjecture is settled only for some special cases. Even a partial solution\footnote{E.\ g. for almost disjoint hypergraphs.  A hypergraph $(V,E)$ is \textit{disjoint} if $A,B\in E\Rightarrow |A\cap B|\leq 1$ if $A\neq B$.}
would be desirable, since one could use Chooser-Picker games as bounds
for what is known as $\alpha-\beta$ pruning of Maker-Breaker games~\cite{CsD}.

Even more importantly, Beck realized that the outcome of a Chooser-Picker game coincides with the outcome of a Maker-Breaker game for some hypergraphs. This correspondence turns out to be extremely fruitful, see \cite{Beck5}. The most striking example is the clique game, where $V(\mathcal H)=E(K_n)$, and $A\in E(\mathcal H)$ iff $A$ is the collection of the edges of $q$-element clique of $K_n$.

Since the Maker-Breaker (and the Maker-Maker) games are PSPACE-complete (see \cite{Schaefer}) it would support both Conjecture~\ref{Bc}, and the above coincidence with Chooser-Picker games to see that the Chooser-Picker or Picker-Chooser games are not easy as well. To prove PSPACE-completeness for positional games is more or less standard, see \cite{Reisch, Byskov}. Here we can prove something weaker because of the asymmetric nature of these games.

\begin{theorem} \label{PCNP}
It is NP-hard to decide the winner in a Picker-Chooser game.
\end{theorem}

\begin{theorem} \label{CPNP}
It is NP-hard to decide the winner in a Chooser-Picker game.
\end{theorem}

In Section~\ref{genpair} we generalize the pairing strategies first formalized by Hales and Jewett \cite{H-J}. As an application, we show there is no pairing strategy for the game ``Snaky,'' see \cite{HararySnaky,HarborthSnaky,SiebenSnaky}. Finally, we compare the actual complexity of these games on a specific hypergraph, the $4 \times 4$ torus in Section~\ref{tori}.

\section{Proofs of Theorems~\ref{PCNP} and \ref{CPNP}}

Both proofs are based on the usual reduction method. We reduce $3-{\rm SAT}$ to Chooser-Picker or Picker-Chooser games.

\medskip

\noindent {\bf Proof of Theorem~\ref{PCNP}.}
Consider an arbitrary CNF formula $\phi(x_1, \dots, x_n) \in {\rm 3-SAT}$. We denote $\phi=C_1\wedge \dots \wedge C_k$, where $C_i=\ell_{i_1}\vee\ell_{i_2}\vee\ell_{i_3}$ and $\ell_{i_j}$ is a literal for $i\in \{1,\ldots,k\}$ and $j=1,2,3$. With a slight abuse of notation, we use $C_i$ also to denote the set of literals in it. That is, if there exists a clause $C_i=x_2 \vee {\bar x_5} \vee x_6$, then we also denote the set $C_i=\{x_2, {\bar x_5}, x_6\}$.

We will exhibit a hypergraph ${\mathcal H}_{\phi}=(V, E)$ such that the Picker-Chooser game is a win for Chooser if and only if $\phi$ is satisfiable.

The vertex set will be $V=\{x_1,\ldots,x_n,{\bar x_1},\ldots,{\bar x_n}\}$. Let $\B\subset 2^V$ have the property that $B\in\B$ if, for all $i\in \{1,\ldots,n\}$, $B$ contains either $x_i$ or ${\bar x_i}$ but not both. The edge set $E$ consists of the sets $A$ such that $A=C_i\cup B$ for some $i$ and some $B\in\B$.

Note that $\B$, and consequently $E$, has a short (polynomial in $\phi$) description even though $|E|\geq |\B|=2^n$.

Claim~1 allows us to restrict our attention to games in which Picker has a specific kind of strategy.
\medskip

\noindent {\bf Claim~1:} If Picker fails to select pairs of the form $\{x_i, {\bar x_i} \}$ in each round, then Chooser has a winning strategy.

\noindent {\bf Proof:} We assume to the contrary: Let $\{x,y\}$ be the first pair selected by Picker such that $\{x,y\}\neq\{x_i,{\bar x_i}\}$ for any $i\in\{1,\ldots,n\}$. In that case, Chooser keeps, say, $x$, and waits until Picker offers up ${\bar x}$ in a pair. In that round, Chooser takes ${\bar x}$, and wins the game, since Picker cannot take any $B\in\B$. This proves Claim~1.\bizveg
\medskip

First we show that if Picker-Chooser on ${\mathcal H}_{\phi}$ is a win for Chooser, then $\phi$ is satisfiable. According to Claim~1, we may assume that Picker's strategy is to select pairs of the form $\{x_i,{\bar x_i}\}$ resulting in the fact that such pairs are shared among Picker and Chooser for all $i$. Assume that Chooser wins the game on ${\mathcal H}_{\phi}$, and set $\hat x_i=1$ if Chooser holds $x_i$, and $\hat x_i=0$ otherwise. Picker holds all elements of some $B\in\B$, so the assumption means that Chooser has an element in each of the $C_i$'s.  That is, $\phi({\hat x_1},\ldots,{\hat x_n})=1$.

Next we show that if $\phi$ is satisfiable, then Picker-Chooser on ${\mathcal H}_{\phi}$ is a win for Chooser.  Since $\phi$ is satisfiable, there exist ${\hat x_1},\ldots,{\hat x_n}$, such that $\phi({\hat x_1},\ldots,{\hat x_n})=1$. Consider the Picker-Chooser game on ${\mathcal H}_{\phi}$. By Claim~1, we may assume that, in each round, Picker offers a pair of the form $\{x_i, {\bar x_i}\}$.  In that case, Chooser takes $x_i$ if and only if ${\hat x_i=1}$, and wins the game. This proves Theorem~\ref{PCNP}.\bizveg

\medskip

\noindent {\bf Proof of Theorem~\ref{CPNP}.}
Let us use the same set-up and notation for the CNF formula $\phi$ as in the proof of Theorem~\ref{PCNP}. We want to define a hypergraph ${\mathcal H}_{\phi}=(V,E)$ such that the Chooser-Picker game on ${\mathcal H}_{\phi}=(V,E)$ is a Picker's win if and only if $\phi$ is satisfiable.

Let the vertex set be $V=\{a_i, b_i, c_i, d_i\}_{i=1}^n$. The edge set, $E$, consists of all edges $A$ such that
\begin{itemize}
   \item $A \subset \{a_i, b_i, c_i, d_i\}$ and $|A|=3$ for some $i \in \{1,\ldots,n\}$,
   \item $A=\{a_i, a_j, a_k, b_i, b_j, b_k \}$ for a clause $C=x_i \vee x_j \vee x_k$,
   \item $A=\{a_i, a_j, a_k, b_i, b_j, c_k \}$ for a clause $C=x_i \vee x_j \vee {\bar x_k}$,
   \item $A=\{a_i, a_j, a_k, b_i, c_j, c_k \}$ for a clause $C=x_i \vee {\bar x_j} \vee {\bar x_k}$,
   \item $A=\{a_i, a_j, a_k, c_i, c_j, c_k \}$ for a clause $C={\bar x_i} \vee {\bar x_j} \vee {\bar x_k}$.
\end{itemize}

Claim~2 allows us to restrict our attention to games in which Chooser has a specific kind of strategy.
\medskip

\noindent {\bf Claim~2:}~\\ \vspace{-0.5cm}
\begin{itemize}
   \item If Picker picks a pair $(x,y)$ such that $\{x,y\}\not\subset\{a_i,b_i,c_i,d_i\}$ for some $i\in\{1,\ldots,n\}$, then Chooser has a winning strategy.
   \item Chooser has an optimal strategy that results in always choosing $a_i$ and always giving $d_i$ to Picker.
\end{itemize}
In particular, this means that we may assume that for all $i$, Picker either picks $\{(a_i,b_i),(c_i,d_i)\}$ or $\{(a_i,c_i),(b_i,d_i)\}$.  Moreover, Chooser will get $a_i$ and Picker will get $d_i$ and each player will get exactly one of $(b_i,c_i)$.

\smallskip

\noindent {\bf Proof:} Suppose Picker offers a pair $(x,y)$ for which $x \in \{a_i, b_i, c_i, d_i\}$ but $y \not \in \{a_i, b_i, c_i, d_i\}$.  Consider the first such instance. In that case, Chooser chooses $x$, and ultimately wins by choosing at least two more elements from $\{a_i,b_i,c_i, d_i\}\setminus\{x\}$, giving Chooser every element of some $A$ of size 3.  So, for all $i$, Picker will pick either $\{(a_i,d_i),(b_i,c_i)\}$ or $\{(a_i,b_i),(c_i,d_i)\}$ or $\{(a_i,c_i),(b_i,d_i)\}$.  Hence, Chooser and Picker will have at least one member of each set of size 3.

However, no $d_i$ appears in any of the sets of size 6 and so if Chooser wins by choosing $d_i$, then he must also win by not choosing $d_i$.  Finally, suppose Picker picks the pair $(a_i,b_i)$ or $(a_i,c_i)$.  Chooser will choose $a_i$ in either case because every $A$ of size 6 that contains either $b_i$ or $c_i$ will also contain $a_i$. So, once again, Chooser can only benefit by choosing $a_i$ over $b_i$ or $c_i$.  Summarizing, if Picker plays optimally; i.e., always taking pairs with the same subscript, then for every winning strategy in which Chooser chooses $d_i$, there exists a winning strategy in which he does not and for every winning strategy in which Chooser does not choose $a_i$, there exists a winning strategy in which he does.

So, we may assume that Picker picks either $\{(a_i,b_i),(c_i,d_i)\}$ or $\{(a_i,c_i),(b_i,d_i)\}$ for all $i$ because if Picker picks $\{(a_i,d_i),(b_i,c_i)\}$, then the outcome is the same except that he cannot control which of $\{b_i,c_i\}$ he will be given by Chooser. This proves Claim~2.\bizveg

\medskip

Now let Picker's $\{(a_i, b_i), (c_i, d_i)\}$ or $\{(a_i, c_i), (b_i, d_i)\}$ moves correspond to setting the value of $x_i=1$ or $x_i=0$, respectively.

First we show that if Chooser-Picker on ${\mathcal H}_{\phi}$ is a win for Picker, then $\phi$ is satisfiable.  We may assume that Chooser plays according to the restrictions imposed by Claim~2. At the end of the game, Picker has exactly one of $\{b_i,c_i\}$.  Chooser has $a_i$ for all $i\in\{1,\ldots,n\}$.  Let ${\hat x_i}=1$ if Picker has $b_i$ and ${\hat x_i}=0$ otherwise. By the construction of ${\mathcal H}_{\phi}$, this means that $\phi({\hat x_1},\ldots,{\hat x_n})=1$.

Next we show that if $\phi$ is satisfiable, then Picker-Chooser on ${\mathcal H}_{\phi}$ is a win for Picker.  Suppose that there is some assignment that $\phi=({\hat x_1},\ldots,{\hat x_n})$. Picker makes sure to get $b_i$ (i.e., Picker picks $\{(a_i,b_i),(c_i,d_i)\}$) if ${\hat x_i}=1$, and makes sure to get $c_i$ (i.e., Picker picks $\{(a_i,c_i),(b_i,d_i)\}$) if ${\hat x_i}=0$.  Because of Claim~2, we may assume that Chooser will always choose $a_i$ for all $i\in\{1,\ldots,n\}$. As a result, Picker will get at least one element from every $A \in E$, and wins the game.  This proves Theorem~\ref{CPNP}.
\bizveg

\medskip

Note that this theorem implies that Chooser-Picker games are NP-hard, even in the case of hypergraphs $(V, E)$, for which $|A| \leq 6$ for all $A \in E$.

\section{Pairing strategies revisited} \label{genpair}

\subsection{Pairing strategies in general}
Pairing strategies appear in a plethora of games, see \cite{B-C-G}. Certain kind of pairing strategies were introduced to the theory of Positional Games by Hales and Jewett in \cite{H-J}. Based on these pairing strategies they proved the following theorem.

\begin{theorem} \cite{H-J} \label{HJ} Breaker wins a Maker-Breaker game on the hypergraph $(V, E)$ if $|\cup_{A \in \G}A| \geq 2|\G|$ for all $\G \subset E$.
\end{theorem}

The idea is to use the celebrated K\H{o}nig-Hall theorem\footnote{A generalized form of this theorem will be spelled out in the next paragraph as Theorem~\ref{infinite}.}, and exhibit a ``double'' {\em system of distinct representatives (SDR)}, in the hypergraph $(V, E)$. A set $X \subset V$ is an SDR if $|X|=|E|$, and there is a bijection $\phi: X \rightarrow E$ such that for all $x \in X$, $x \in \phi(x)$. If $X$ and $Y$ are SDR's of $(V, E)$ with the bijections $\phi$ and $\psi$ where $X \cap Y=\emptyset$, then $\rho=\psi^{-1}\phi$ is a bijection $\rho: X \rightarrow Y$. Breaker, even as a second player, wins by using $\rho$. That is, Breaker takes $\rho(x)$ [takes $\rho^{-1}(y)$] if Maker takes an $x \in X$ [a $y \in Y$], and an arbitrary untaken element $v \in V$ if Maker takes a $w \in V\setminus (X \cup Y)$.

While Theorem~\ref{HJ} works fine for some games, it has its drawbacks. It rarely gives sharp results, which is not surprising considering the PSPACE-completeness of those games. Another problem is that the K\H{o}nig-Hall theorem (and consequently Theorem~\ref{HJ}) applies only to finite hypergraphs. A remedy for this is a lesser known theorem of Marshall Hall Jr., that requires only the local finiteness of the hypergraph $(V, E)$. We say that $(V, E)$ is {\em locally finite} if ${\rm deg}(x):=|\{A: x \in A \in E\}| < \infty$ for all $x \in V$.

\begin{theorem} \cite{MHall} \label{infinite} There is a SDR in a locally finite hypergraph $(V, E)$ iff $|\cup_{A \in \G}A| \geq |\G|$ for all $\G \subset E$.
\end{theorem}

Still, Theorem~\ref{HJ} does not apply directly if $|V| < 2|E|$, for instance, one must use other ideas to tackle the $k$-in-a-row games in two or in higher dimensions, see \cite{Plu2}.

\begin{definition} The bijection $\rho: X \rightarrow Y$, where $X \cap Y=\emptyset$ and $X, Y \subset V$, is a {\bf winning pairing strategy} for Breaker in the Maker-Breaker game on hypergraph $(V, E)$ if for all $A \in E$ there is an $x \in X$ such that $\{x, \rho(x)\} \subset A$.
\end{definition}

Of course, we assume that both the function $\rho$ and the decision problem that determining whether any set $Y \subset V$ has the property that $Y \subset A \in E$ are computable in polynomial time in the size of description of $(V, E)$. (For the sake of simplicity we consider only the case when both $V$ and $E$ are finite.) Having the bijection $\rho$, Breaker wins by taking $\rho(x)$ [taking $\rho^{-1}(y)$] if Maker's last move was $x \in X$ [was $y \in Y$]. To decide the existence of $\rho$ is not easy in general. Let us denote the class of hypergraphs for which Breaker has a winning pairing strategy by $\B$.

\begin{theorem} \label{general}
Determining whether a hypergraph is in $\B$ is NP-complete.
\end{theorem}

\noindent{\bf Proof.} Given a bijection $\rho$ that witnesses a winning pairing strategy, one checks for an $A \in E$ if there is an $x \in X$ such that $\{x, \rho(x)\} \subset A$. For any pair $(A,x)$ it can be done in polynomial time, and $|E||V|$ is an upper bound on the number of such pairs. Consequently, $\B \in {\rm NP}$.

To show that $\B$ is NP-hard one can use basically the same argument as in the proof of Theorem~\ref{CPNP}. There is, however, a simpler reduction. Let $\phi$ be an arbitrary CNF in $3$-SAT. We construct a hypergraph $\mathcal{H}_{\phi}=(V,E)$ such that $V=\{r_i, b_i, p_i\}_{i=1}^n$ and the edge set, $E$, consists of all edges $A$ such that
\begin{itemize}
   \item $A$ is $ \{r_i, b_i, p_i\}$ for all $i\in\{1,\ldots,n\}$,
   \item $A=\{p_i, r_i, p_j, r_j, p_k, r_k\}$ for a clause $C=x_i \vee x_j \vee x_k$,
   \item $A=\{p_i, r_i, p_j, r_j, p_k, b_k\}$ for a clause $C=x_i \vee x_j \vee {\bar x_k}$,
   \item $A=\{p_i, r_i, p_j, b_j, p_k, b_k\}$ for a clause $C=x_i \vee {\bar x_j} \vee {\bar x_k}$,
   \item $A=\{p_i, b_i, p_j, b_j, p_k, b_k\}$ for a clause $C={\bar x_i} \vee {\bar x_j} \vee {\bar x_k}$.
\end{itemize}

A winning pairing strategy for Breaker cannot contain both $\{p_i, r_i\}$ or $\{p_i, b_i\}$ for $1\leq i\leq n$, because the strategy is a bijection.  But such a strategy must contain one of $\{p_i, r_i\}$ or $\{p_i, b_i\}$ in order to have at least one pair of the form $\{x,\rho(x)\}$ in each of the edges of size 3. Let $x_i=1$ if $\{p_i, r_i\}$ is present, while $x_i=0$ otherwise. As a result, a clause $C$ associated to its corresponding set $A$ of size $6$ is satisfied if and only if $A$ contains a pair. \qed

\medskip

{\noindent \bf Remarks.} If the hypergraph $(V, E)$ is almost disjoint, then Breaker has a winning pairing strategy iff $|\cup_{A \in \G}A| \geq 2|\G|$ for all $\G \subset E$, that is one gets back the assumption of Theorem~\ref{HJ}. This case can be decided in polynomial time in the description of $(V, E)$. As in Theorem~\ref{CPNP}, $\B$ is NP-complete for hypergraphs $(V, E)$, where $|A| \leq 6$ for $A \in E$. A result of Hegyh\'ati \cite{Hegyhati} implies that the existence of a winning pairing strategy can be decided in polynomial time if $|A| \leq 3$ for $A \in E$. The cases when $|A| \leq 4$ or $|A| \leq 5$, to the best of our knowledge, are open.

\subsection{Applications for k-in-a-row and Snaky}

Let $d_2$ be the maximum pair degree in $(V, E)$, that is $d_2=\max_{x \neq y}d_2(x, y)$, where $d_2(x, y)=|\{A: \{x, y\} \subset A \in E\}|$.

\begin{proposition} \label{trivial}
If Breaker has a winning pairing strategy then $d_2|X|/2 \geq |\G|$ must hold for all $X \subset V$, where $\G=\{A: A \in E, A \subset X\}$.
\end{proposition}

\noindent{\bf Proof.} Simply locate the pairs in the winning pairing strategy.  There are at most $|X|/2$ such pairs, which are disjoint.  Each pair will be a subset of at most $d_2$ edges.  Since each edge of $\G$ must have a pair as a subset, the number of edges must be at most $d_2|X|/2$. \qed
\medskip

Now we can explain why pairing strategies can work for the game $k$-in-a-row for sufficiently large $n$ only if $k \geq 9$, see \cite{B-C-G}. In the $k$-in-a-row game, $d_2=k-1$, and if $X$ is an $n \times n$ board, then $|\G|=4n^2+O(kn)$. By Proposition~\ref{trivial}, we have $(k-1)n^2/2 \geq 4n^2+O(kn)$; that is, $k \geq 9+o(n)$.
\smallskip

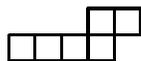
\begin{figure}[htbp]
\begin{picture}(60,30)(0,0)
\thicklines
\put(0,0){\line(1,0){40}}
\put(0,10){\line(1,0){50}}
\put(30,20){\line(1,0){20}}
\multiput(0,0)(10,0){5}{\put(0,0){\line(0,1){10}}}
\multiput(30,10)(10,0){3}{\put(0,0){\line(0,1){10}}}
\end{picture}
\caption{The polyomino Snaky. The ``head'' is the pair of cells in the upper row.  The ``body'' is the set of four consecutive cells in the lower row.}
\label{snaky}
\end{figure}

Another example in which we can use this ideas is the polyomino game Snaky, which were examined by Harary~\cite{HararySnaky}, Harborth and Seeman~\cite{HarborthSnaky}, and Sieben~\cite{SiebenSnaky}. This game is a Maker-Breaker game in which the board consists of the cells of the infinite grid and Maker's goal is to occupy all of the cells in an isomorphic copy of the polyomino Snaky, shown in Figure~\ref{snaky}.

Using a computer search, Harborth and Seeman \cite{HarborthSnaky} showed that there is no pairing strategy for Breaker in this game. We give a computer-free proof for their statement:

\begin{theorem} \cite{HarborthSnaky} \label{nomatch}
Breaker has no pairing strategy to avoid the isomorphic copies of the polyomino ``Snaky.''
\end{theorem}

\noindent{\bf Proof.} Assume to the contrary that there is a winning pairing $\rho$ for Breaker. Let $P_\ell$ be the polyomino which consists of $\ell$ consecutive squares of the table.

First we show that $\rho$ cannot be a pairing for the polyomino $P_4$. Let us assume that $\rho$ is such a pairing, and consider an $n \times n$ board $X$ such that the edges of $\G$ consist of the $P_4$'s on $X$. Since $d_2=3$, Proposition~\ref{trivial} gives that $3n^2/2 \geq 2n^2+O(n)$, which is a contradiction if $n$ is sufficiently large.

On the other hand, if $\rho$ is a pairing for Snaky, then we will show that it must be a pairing for $P_5$. To see this, we assign labels to the cells such that cells receive the same label iff they are paired by $\rho$. Let us take the longest set of consecutive cells $R$ such that no labels are repeated on $R$. We may assume that either those labels are $1, \dots, \ell$ for some $\ell \geq 5$, or $R$ is infinite.

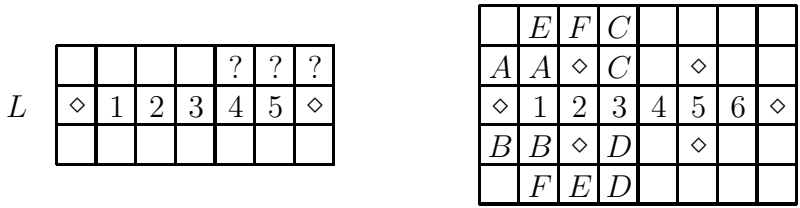
\begin{figure}[htbp]
\begin{picture}(300,90)(0,0)
\thicklines
\multiput(0,0)(0,15){4}{\put(0,15){\line(1,0){105}}}
\multiput(0,0)(15,0){8}{\put(0,15){\line(0,1){45}}}
\put(-20, 37){\makebox(10, 0){$L$}}
\put(18, 37){\makebox(10, 0){$1$}}
\put(33, 37){\makebox(10, 0){$2$}}
\put(48, 37){\makebox(10, 0){$3$}}
\put(63, 37){\makebox(10, 0){$4$}}
\put(78, 37){\makebox(10, 0){$5$}}
\put(63, 52){\makebox(10, 0){$?$}}
\put(78, 52){\makebox(10, 0){$?$}}
\put(93, 52){\makebox(10, 0){$?$}}
\put(3, 37){\makebox(10, 0){$\diamond$}}
\put(93, 37){\makebox(10, 0){$\diamond$}}
\multiput(160,0)(0,15){6}{\put(0,0){\line(1,0){120}}}
\multiput(160,0)(15,0){9}{\put(0,0){\line(0,1){75}}}
\put(178, 37){\makebox(10, 0){$1$}}
\put(193, 37){\makebox(10, 0){$2$}}
\put(208, 37){\makebox(10, 0){$3$}}
\put(223, 37){\makebox(10, 0){$4$}}
\put(238, 37){\makebox(10, 0){$5$}}
\put(253, 37){\makebox(10, 0){$6$}}
\put(268, 37){\makebox(10, 0){$\diamond$}}
\put(163, 52){\makebox(10, 0){$A$}}
\put(178, 52){\makebox(10, 0){$A$}}
\put(163, 22){\makebox(10, 0){$B$}}
\put(178, 22){\makebox(10, 0){$B$}}
\put(208, 52){\makebox(10, 0){$C$}}
\put(208, 22){\makebox(10, 0){$D$}}
\put(163, 37){\makebox(10, 0){$\diamond$}}
\put(193, 52){\makebox(10, 0){$\diamond$}}
\put(238, 52){\makebox(10, 0){$\diamond$}}
\put(193, 22){\makebox(10, 0){$\diamond$}}
\put(238, 22){\makebox(10, 0){$\diamond$}}
\put(178, 67){\makebox(10, 0){$E$}}
\put(208, 67){\makebox(10, 0){$C$}}
\put(193, 7){\makebox(10, 0){$E$}}
\put(208, 7){\makebox(10, 0){$D$}}
\put(193, 67){\makebox(10, 0){$F$}}
\put(178, 7){\makebox(10, 0){$F$}}
%
\end{picture}
\caption{The cases $\ell=5$ and $\ell=6$.}
\label{five}
\end{figure}

We first consider the case $\ell=5$, and in doing so let us refer to a cell of the grid by its lower left lattice point. If $\rho$ is not a pairing for $P_5$, then we may assume, without loss of generality, that the set of cells $L=\{(1,0),\ldots,(5,0)\}$ contains no pairs. These cells are labeled by $1, \dots, 5$ on the left-hand side of Figure~\ref{five}. Since $\ell=5$, the both the cells $(0, 0)$ and $(6, 0)$ are in a pair with some cell of $L$. (We indicate the cells that have indices which matching with an element of $L$ by a diamond, otherwise by capital letters.) This leaves only three elements of $L$ that can be matched with a cell the rows above and below of $L$.

Consider the Snakys that have four cells in $L$. The head of the snake will have two cells in one of 4 disjoint sets of three consecutive cells in the row above or the row below $L$. Without loss of generality, we may assume that the three consecutive cells $\{(4,1),(5,1),(6,1)\}$.  That is, no cell of $L$ is matched by the cells $\{(4,1),(5,1),(6,1)\}$, labeled by ``?'' in Figure~\ref{five}. But in that case $\rho$ should contain, as pairs, both $\{(4,1),(5,1)\}$ and $\{(5,1),(6,1)\}$, which is impossible. So we may assume that $\ell>5$.~\\

\noindent\textbf{Remark.} In the case that $\ell>5$, or $\ell$ is infinite, we again have a set $L$ containing no pairs such that $|L|=\ell$. Every three consecutive cells in the rows above and below $L$ must contain at least one cell whose label is matched to a cell of $L$, otherwise we finish the argument as in case $\ell=5$. Here by ``the rows above and below $L$'' we mean sets that extend one cell longer than the end of $L$ if $L$ is finite or if $L$ terminates in one direction.~\\

Second is the case of $\ell = 6$ and we may assume that $\{(1,0),\ldots,(6,0)\}$ receive distinct labels. We will show that the only possible pattern is shown in the right-hand side of Figure~\ref{five}.  There are diamonds in the cells $(0,0)$ and $(7,0)$.  Four diamonds remain to be placed and each set of three consecutive cells above and below $L$.  The only possible locations do to so are $(2,\pm 1)$ and $(5,\pm 1)$. This ensures that $\{(0,1),(1,1)\}$ and $\{(0,-1),(1,-1)\}$ form pairs, which we label with ``$A$'' and ``$B$'', respectively.

Note that neither diamonds above and below the cell ``2'' can also be labeled by ``2'', otherwise the diamond, its right neighbor, and the cells $3,4,5,6$ would be a pairing-free Snaky. The cells above and below the cell ``3'' are labeled ``$C$'' and ``$D$'', respectively. At this moment $C$ could be equal to $D$. However, if we consider a standing Snaky on the cells $\{(1,2),(1,1),(2,1),(2,0),(2,-1),(2,-2)\}$, the only unpaired cells are those that are labeled with ``$E$''. If we consider a standing Snaky with the same body and the head towards the upper right, the only unpaired cells are those labeled ``$C$'' in the right-hand side of Figure~\ref{five}.  Symmetrically, we may assign labels ``$D$'' and ``$F$'' as shown in the figure. This, however, leads to a contradiction, since there would be a pairing-free Snaky again. In particular, the upper $E$ and $F$ cells make the head, and the body consists of the diamond above the cell ``2'', the cell of the lower $C$, the empty cell above ``4'' and the diamond above the cell ``5''. So, we may assume that $\ell>6$.~\\

The third case, where $\ell=7$, is impossible since the rows above and below $L$ should contain three diamonds each to avoid the snakes and two are needed to the right and left of $L$.  This totals at least 8, more than the $7$ that are available.~\\

\begin{figure}[htbp]
\begin{picture}(300,90)(0,0)
\thicklines
\multiput(-30,0)(0,15){4}{\put(0,15){\line(1,0){150}}}
\multiput(-30,0)(15,0){11}{\put(0,15){\line(0,1){45}}}
\put(-50, 37){\makebox(10, 0){$L$}}
\put(-12, 37){\makebox(10, 0){$1$}}
\put(3, 37){\makebox(10, 0){$2$}}
\put(18, 37){\makebox(10, 0){$3$}}
\put(33, 37){\makebox(10, 0){$4$}}
\put(48, 37){\makebox(10, 0){$5$}}
\put(63, 37){\makebox(10, 0){$6$}}
\put(78, 37){\makebox(10, 0){$7$}}
\put(93, 37){\makebox(10, 0){$8$}}
\put(-27, 37){\makebox(10, 0){$\diamond$}}
\put(108, 37){\makebox(10, 0){$\diamond$}}
\put(18, 52){\makebox(10, 0){$\diamond$}}
\put(63, 52){\makebox(10, 0){$\diamond$}}
\put(33, 52){\makebox(10, 0){$A$}}
\put(48, 52){\makebox(10, 0){$A$}}
\multiput(160,0)(0,15){4}{\put(0,15){\line(1,0){165}}}
\multiput(160,0)(15,0){12}{\put(0,15){\line(0,1){45}}}
\put(140, 37){\makebox(10, 0){$L$}}
\put(178, 37){\makebox(10, 0){$1$}}
\put(193, 37){\makebox(10, 0){$2$}}
\put(208, 37){\makebox(10, 0){$3$}}
\put(223, 37){\makebox(10, 0){$4$}}
\put(238, 37){\makebox(10, 0){$5$}}
\put(253, 37){\makebox(10, 0){$6$}}
\put(268, 37){\makebox(10, 0){$7$}}
\put(283, 37){\makebox(10, 0){$8$}}
\put(298, 37){\makebox(10, 0){$9$}}
\put(223, 52){\makebox(10, 0){$\diamond$}}
\put(238, 52){\makebox(10, 0){$\diamond$}}
\put(253, 52){\makebox(10, 0){$\diamond$}}
\put(223, 22){\makebox(10, 0){$\diamond$}}
\put(238, 22){\makebox(10, 0){$\diamond$}}
\put(253, 22){\makebox(10, 0){$\diamond$}}
\end{picture}
\caption{The cases $\ell=8$ and $\ell\geq 9$.}
\label{eight}
\end{figure}
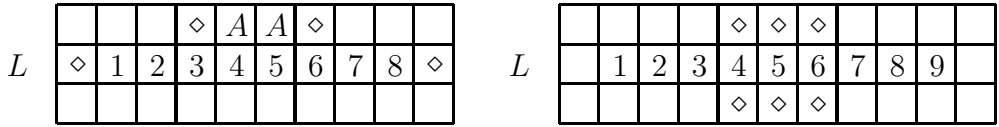

In the fourth case, where $\ell=8$, we have at most eight diamonds around $L$, two of those at the ends, and every three consecutive cells above and below $L$ containing at least one diamond. So, there are ten cells above $L$ and ten cells below $L$ to receive the remaining 6 diamonds.  There must be one in the three leftmost cells above $L$, in the three rightmost cells above $L$, in the three leftmost cells below $L$ and in the three rightmost cells below $L$.  Only two diamonds remain.  One must be above one of the cells labeled ``$3$'', ``$4$'', ``$5$'' or ``$6$''. A diamond cannot be above the cell labeled ``$4$'' or ``$5$'' because for the two Snakys with heads equal to $\{(4,1),(5,1)\}$ and bodies in $L$, the diamond either represents one of $\{1,2,3,4\}$ or one of $\{5,6,7,8\}$.  Hence, one of these Snakys must be pairing-free. As a result, the cells $\{(4,1),(5,1)\}$ must be paired with each other and so we label them with ``$A$''. See the diagram in the left-hand side of Figure~\ref{eight}.  Because every three consecutive cells must contain at least one diamond, the cells above the cells labeled ``3'' and ``6'' are labeled with a diamond.  This is a contradiction to the fact that only one diamond can be above these cells. So, we may assume that $\ell>8$.~\\

In the fifth case, where $\ell \geq 9$ and is finite, all cells above and below the cells $4, \dots, \ell -3$, the ``critical region'', must be diamonds. It is the same idea as in the previous case: If, say the cell above ``4'', is $A$, then so is the cell above ``5''. But the same is true for the cells above ``5'' and ``6''. Not only must the cells in the critical region be diamonds, there must be a total of at least 4 more above at below $L$ to cover all of the triples of consecutive cells.  With the additional two on the endpoints, there must be at least $2(\ell-6)+4 + 2$ diamonds, that is impossible, given that the total number of diamonds is at most $\ell$, which is at least 9.~\\

Finally, suppose $L$ is infinite.  Take 13 consecutive cells of $L$, call it $L'$. In the critical region of $L'$ there must be $2(13-6)=14$ cells with diamonds, but they must repeat the labels in the cells of $L'$, a contradiction.  This concludes the proof of the fact that a pairing for Snaky must be a pairing for $P_5$.~\\

\begin{figure}[htbp]
\begin{picture}(300,90)(0,0)
\thicklines
\multiput(0,0)(0,20){5}{\put(0,0){\line(1,0){80}}}
\multiput(0,0)(20,0){5}{\put(0,0){\line(0,1){80}}}
\multiput(10,0)(40,0){2}{\multiput(0,0)(20,20){2}
{\multiput(0,0)(0,40){2}{\put(0,0){\line(0,1){20}}}}}
\multiput(20,10)(40,0){2}{\multiput(0,0)(-20,20){2}
{\multiput(0,0)(0,40){2}{\put(0,0){\line(1,0){20}}}}}
\multiput(210,10)(10,10){7}{\multiput(0,0)(0,10){2}{{\put(0,0){\line(1,0){20}}}}}
\multiput(240,0)(10,10){7}{\multiput(0,0)(0,10){2}{{\put(0,0){\line(1,0){20}}}}}
\multiput(180,20)(10,10){7}{\multiput(0,0)(0,10){2}{{\put(0,0){\line(1,0){20}}}}}
\multiput(200,10)(30,-10){2}{\put(0,0){\line(1,0){10}}}
\multiput(260,90)(30,-10){2}{\put(0,0){\line(1,0){10}}}
\multiput(200,10)(10,10){7}{\multiput(0,0)(10,0){2}{{\put(0,0){\line(0,1){20}}}}}
\multiput(230,0)(10,10){7}{\multiput(0,0)(10,0){2}{{\put(0,0){\line(0,1){20}}}}}
\multiput(180,20)(10,10){7}{{\put(0,0){\line(0,1){10}}}}
\multiput(260,0)(10,10){7}{{\put(0,0){\line(0,1){10}}}}
\end{picture}
\caption{The parings $T_1$ and $T_2$.}
\label{nofive}
\end{figure}
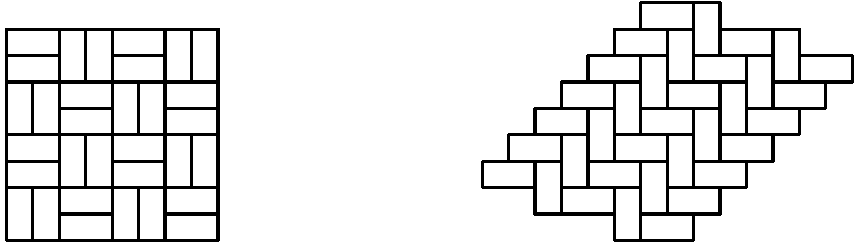

We exhibit two pairings for $P_5$. The pairing $T_1$ is like a chessboard, where the fields are $2 \times 2$, and alternately packed by a standing and lying pairs of dominoes as in the left-hand side of Figure~\ref{nofive}. The pairing $T_2$ is like an infinite zipper, repeated in both directions, see the right-hand side of Figure~\ref{nofive}.

\begin{lemma} \label{tessalation}
A pairing for $P_5$ is either the translated and rotated copy of $T_1$ or $T_2$.
\end{lemma}

\noindent{\bf Proof.} Let us consider a pairing, $\rho$, for $P_5$. A pair $\{x, \rho(x)\}$ is {\em good} if $x$ and $\rho(x)$ are neighboring cells. If $\{x, \rho(x)\}$ is good, then $d_2(x, \rho(x))=4$, otherwise it is smaller. The number of $P_5$'s are $2n^2+O(n)$ on an $n \times n$ sub-board $X$, so Proposition~\ref{trivial} implies that at all but $O(n)$ pairs on $X$ are good. It follows that, if $n$ is sufficiently large, then there is a $Y \subset X$, $k \times k$ square sub-board that contains only good pairs. I.\ e. this $k\times k$ sub-board is paired by dominoes.

There are either two dominoes meeting at their longer sides, or the two long sides meet but are offset by one unit. In these cases the immediate neighboring dominoes are forced to be in the pattern of $T_1$ or $T_2$, respectively.

\begin{figure}[htbp]
\begin{picture}(300,120)(0,-20)
\thicklines
\multiput(-40,0)(0,20){5}{\put(0,0){\line(1,0){80}}}
\multiput(-40,0)(20,0){5}{\put(0,0){\line(0,1){80}}}
\multiput(-10,0)(40,0){2}{\multiput(0,0)(20,20){2}
{\multiput(-20,0)(0,40){2}{\put(0,0){\line(0,1){20}}}}}
\multiput(0,10)(40,0){2}{\multiput(0,0)(-20,20){2}
{\multiput(-20,0)(0,40){2}{\put(0,0){\line(1,0){20}}}}}
\put(-10, 30){\makebox(10, 0){{\tiny 1}}}
\put(-20, 30){\makebox(10, 0){{\tiny 1}}}
\put(5, 35){\makebox(10, 0){{\tiny 2}}}
\put(5, 25){\makebox(10, 0){{\tiny 2}}}
\put(-35, 35){\makebox(10, 0){{\tiny 2}}}
\put(-35, 25){\makebox(10, 0){{\tiny 2}}}
\put(10, 50){\makebox(10, 0){{\tiny 3}}}
\put(0, 50){\makebox(10, 0){{\tiny 3}}}
\put(10, 10){\makebox(10, 0){{\tiny 3}}}
\put(0, 10){\makebox(10, 0){{\tiny 3}}}
\put(-30, 50){\makebox(10, 0){{\tiny 3}}}
\put(-40, 50){\makebox(10, 0){{\tiny 3}}}
\put(-30, 10){\makebox(10, 0){{\tiny 3}}}
\put(-40, 10){\makebox(10, 0){{\tiny 3}}}
\put(-15, 55){\makebox(10, 0){{\tiny 4}}}
\put(-15, 45){\makebox(10, 0){{\tiny 4}}}
\put(-15, 15){\makebox(10, 0){{\tiny 4}}}
\put(-15, 5){\makebox(10, 0){{\tiny 4}}}
\put(25, 55){\makebox(10, 0){{\tiny 4}}}
\put(25, 45){\makebox(10, 0){{\tiny 4}}}
\put(25, 15){\makebox(10, 0){{\tiny 4}}}
\put(25, 5){\makebox(10, 0){{\tiny 4}}}
\put(-10, 70){\makebox(10, 0){{\tiny 5}}}
\put(-20, 70){\makebox(10, 0){{\tiny 5}}}
\put(30, 30){\makebox(10, 0){{\tiny 5}}}
\put(20, 30){\makebox(10, 0){{\tiny 5}}}
\put(30, 70){\makebox(10, 0){{\tiny 5}}}
\put(20, 70){\makebox(10, 0){{\tiny 5}}}
\put(5, 75){\makebox(10, 0){{\tiny 6}}}
\put(5, 65){\makebox(10, 0){{\tiny 6}}}
\put(-35, 75){\makebox(10, 0){{\tiny 6}}}
\put(-35, 65){\makebox(10, 0){{\tiny 6}}}
%
\put(90,100){\line(1,0){10}}
\put(80,60){\line(1,0){20}}
\multiput(80,70)(10,10){3}{\put(0,0){\line(1,0){30}}}
\put(130,100){\line(1,0){10}}
\put(80,20){\line(1,0){20}}
\multiput(80,30)(10,10){7}{\put(0,0){\line(1,0){30}}}
\put(170,100){\line(1,0){10}}
\put(90,-10){\line(1,0){20}}
\put(110,-20){\line(1,0){10}}
\multiput(90,0)(10,10){9}{\put(0,0){\line(1,0){30}}}
\put(180,90){\line(1,0){20}}
\multiput(120,-10)(10,10){7}{\put(0,0){\line(1,0){30}}}
\put(190,60){\line(1,0){20}}
\put(150,-20){\line(1,0){10}}
\multiput(160,-10)(10,10){3}{\put(0,0){\line(1,0){30}}}
\put(190,20){\line(1,0){20}}
\put(190,-20){\line(1,0){10}}
%
%
\put(80,60){\line(0,1){10}}
\put(90,70){\line(0,1){30}}
\put(100,80){\line(0,1){20}}
\multiput(90,30)(10,10){5}{\put(0,0){\line(0,1){30}}}
\put(140,80){\line(0,1){20}}
\put(80,20){\line(0,1){10}}
\multiput(90,-10)(10,10){9}{\put(0,0){\line(0,1){30}}}
\put(110,-20){\line(0,1){20}}
\put(180,80){\line(0,1){20}}
\multiput(120,-20)(10,10){9}{\put(0,0){\line(0,1){30}}}
\put(150,-20){\line(0,1){20}}
\multiput(160,-20)(10,10){5}{\put(0,0){\line(0,1){30}}}
\put(210,50){\line(0,1){10}}
\put(190,-20){\line(0,1){20}}
\put(200,-20){\line(0,1){30}}
\put(210,10){\line(0,1){10}}
%
%
\put(130, 10){\makebox(10, 0){{\tiny 4}}}
\put(170, 10){\makebox(10, 0){{\tiny 8}}}
\put(115, 15){\makebox(10, 0){{\tiny 8}}}
\put(155, 15){\makebox(10, 0){{\tiny 7}}}
\put(140, 20){\makebox(10, 0){{\tiny 2}}}
\put(125, 25){\makebox(10, 0){{\tiny 4}}}
\put(165, 25){\makebox(10, 0){{\tiny 6}}}
\put(110, 30){\makebox(10, 0){{\tiny 9}}}
\put(150, 30){\makebox(10, 0){{\tiny 3}}}
\put(135, 35){\makebox(10, 0){{\tiny 1}}}
\put(175, 35){\makebox(10, 0){{\tiny 6}}}
\put(120, 40){\makebox(10, 0){{\tiny 5}}}
\put(160, 40){\makebox(10, 0){{\tiny 5}}}
\put(105, 45){\makebox(10, 0){{\tiny 6}}}
\put(145, 45){\makebox(10, 0){{\tiny 1}}}
\put(130, 50){\makebox(10, 0){{\tiny 3}}}
\put(170, 50){\makebox(10, 0){{\tiny 9}}}
\put(115, 55){\makebox(10, 0){{\tiny 6}}}
\put(155, 55){\makebox(10, 0){{\tiny 4}}}
\put(140, 60){\makebox(10, 0){{\tiny 2}}}
\put(125, 65){\makebox(10, 0){{\tiny 7}}}
\put(165, 65){\makebox(10, 0){{\tiny 8}}}
\put(110, 70){\makebox(10, 0){{\tiny 8}}}
\put(150, 70){\makebox(10, 0){{\tiny 4}}}
%
\multiput(300,10)(10,0){2}{\multiput(0,0)(0,40){2}{{\put(0,0){\line(0,1){20}}}}}
%
\put(290,30){\line(0,1){10}}
\put(310,30){\line(0,1){10}}
\put(300,40){\line(0,1){10}}
\put(320,40){\line(0,1){10}}
\multiput(290,30)(0,10){2}{\multiput(0,0)(10,10){2}{{\put(0,0){\line(1,0){20}}}}}
\put(300,10){\line(1,0){10}}
\put(300,70){\line(1,0){10}}
\put(295, 35){\makebox(10, 0){{\tiny 1}}}
\put(305, 45){\makebox(10, 0){{\tiny 1}}}
\put(300, 20){\makebox(10, 0){{\tiny 2}}}
\put(300, 60){\makebox(10, 0){{\tiny 2}}}
%
%
%
\thinlines
\put(310,20){\line(1,0){12}}
\put(322,20){\line(0,1){50}}
\put(322,70){\line(-1,0){12}}
\put(300,60){\line(-1,0){12}}
\put(288,60){\line(0,-1){50}}
\put(288,10){\line(1,0){12}}
\end{picture}
\caption{The forcing for pairs and filling.}
\label{forcing}
\end{figure}
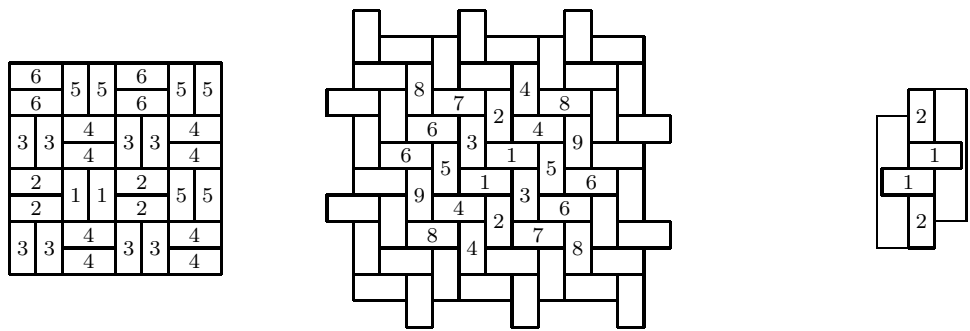

We will show that if we have a large enough pattern of dominoes, then the pairs in the neighboring cells are forced to be in either $T_1$ or $T_2$.
First suppose that, within the pattern tiled by dominoes that two dominoes share a long edge, as in the dominoes labeled with ``1'' in the left-hand side of Figure~\ref{forcing}. Since the pairs can only occur as dominoes, we can use horizontal $P_5$'s to ensure the pairing is oriented as in the dominoes labeled with ``2''.  Vertical $P_5$'s ensure the orientations of the dominoes labeled ``3''.  We can continue in this fashion, getting the $8\times 8$ pattern in the left-hand side of Figure~\ref{forcing}.  Once this is determined, one can extend the pattern to a larger rectangle, forcing not just the domino condition, but the $T_1$ pattern itself.  This can be seen by first taking horizontal $P_5$'s in rows 1,2,5,6 that have two cells outside of the pattern.  Then taking vertical $P_5$'s in columns 9,10, the pattern can be extended to an $8\times 10$ rectangle.  This can be continued ad infinitum, showing that the entire $n\times n$ board must be in the pattern $T_1$.

Next, suppose that whenever two dominoes meet at their long edge in the sub-board, that they are offset by one unit, since two dominoes meeting at their long edge will force the pattern $T_1$.  See dominoes labeled ``1'' in the diagrams in the center or the right-hand side of Figure~\ref{forcing}.  The pairs must occur as dominoes and so vertical $P_5$'s ensure that the dominoes labeled with ``2'' are placed in that location. Now, consider the right-hand side of Figure~\ref{forcing}.  Two $P_5$'s are indicated by thin lines.  Since the dominoes cannot share a long side, this forces the placement of the dominoes labeled with ``3''.

In fact, if we know that a sub-board is tiled with dominoes that do not share a long edge, then the configuration must be that of $T_2$. It remains to show that if we have a large enough fragment of $T_2$ in a sub-board, then, even if the board is not guaranteed to be tiled with dominoes, it must be completed to a $T_2$ pattern. The other pairs are forced even without the assumption that those are in dominoes, since the otherwise a $P_5$ containing no pair would arise.

To see how we can use this sub-board to extend $T_2$ to the whole board, we first show in the center of Figure~\ref{forcing} how enough pairs can be formed under the assumption that every pair forms a domino and no pair of dominoes can share a long edge.  The numbers show the order in which dominoes can be taken.  Then, in Figure~\ref{77to99} we show how, under no assumptions that the pairs occur as dominoes, that the dominoes that cover the $7\times 7$ board can be extended to cover a $9\times 9$ board.  Again, the numbers show the order in which dominoes can be taken.

The general approach is that one can force new horizontal dominoes in every third row that touch the left and right border of the small square and vertical dominoes in every third column that touch the top and bottom border.  From there, the rest of the larger square is easy to complete.  This can continue \textit{ad infinitum} until the board is filled. This concludes the proof of Lemma~\ref{tessalation}.\footnote{How about this? --RM} \qed

\newcommand{\boxfill}{$\cdot$}
\begin{figure}[htbp]
\begin{picture}(130,120)(80,-20)
\thicklines
%
\put(90,100){\line(1,0){10}}
\put(80,60){\line(1,0){20}}
\multiput(80,70)(10,10){3}{\put(0,0){\line(1,0){30}}}
\put(130,100){\line(1,0){10}}
\put(80,20){\line(1,0){20}}
\multiput(80,30)(10,10){7}{\put(0,0){\line(1,0){30}}}
\put(170,100){\line(1,0){10}}
\put(90,-10){\line(1,0){20}}
\put(110,-20){\line(1,0){10}}
\multiput(90,0)(10,10){9}{\put(0,0){\line(1,0){30}}}
\put(180,90){\line(1,0){20}}
\multiput(120,-10)(10,10){7}{\put(0,0){\line(1,0){30}}}
\put(190,60){\line(1,0){20}}
\put(150,-20){\line(1,0){10}}
\multiput(160,-10)(10,10){3}{\put(0,0){\line(1,0){30}}}
\put(190,20){\line(1,0){20}}
\put(190,-20){\line(1,0){10}}
%
%
\put(80,60){\line(0,1){10}}
\put(90,70){\line(0,1){30}}
\put(100,80){\line(0,1){20}}
\multiput(90,30)(10,10){5}{\put(0,0){\line(0,1){30}}}
\put(140,80){\line(0,1){20}}
\put(80,20){\line(0,1){10}}
\multiput(90,-10)(10,10){9}{\put(0,0){\line(0,1){30}}}
\put(110,-20){\line(0,1){20}}
\put(180,80){\line(0,1){20}}
\multiput(120,-20)(10,10){9}{\put(0,0){\line(0,1){30}}}
\put(150,-20){\line(0,1){20}}
\multiput(160,-20)(10,10){5}{\put(0,0){\line(0,1){30}}}
\put(210,50){\line(0,1){10}}
\multiput(190,-10)(10,10){3}{\put(0,0){\line(0,1){10}}}
\put(190,-20){\line(0,1){20}}
\put(200,-20){\line(0,1){30}}
\put(210,10){\line(0,1){10}}
%
%
%
\put(120, 0){\makebox(10, 0){{\tiny 10}}}
\put(160, 0){\makebox(10, 0){{\tiny 10}}}
\put(105, 5){\makebox(10, 0){{\tiny 12}}}
\put(145, 5){\makebox(10, 0){{\tiny 11}}}
\put(185, 5){\makebox(10, 0){{\tiny 12}}}
\put(100, 20){\makebox(10, 0){{\tiny 13}}}
\put(180, 20){\makebox(10, 0){{\tiny 13}}}
\put(95, 35){\makebox(10, 0){{\tiny 10}}}
\put(185, 45){\makebox(10, 0){{\tiny 10}}}
\put(100, 60){\makebox(10, 0){{\tiny 13}}}
\put(180, 60){\makebox(10, 0){{\tiny 13}}}
\put(95, 75){\makebox(10, 0){{\tiny 12}}}
\put(135, 75){\makebox(10, 0){{\tiny 11}}}
\put(175, 75){\makebox(10, 0){{\tiny 12}}}
\put(120, 80){\makebox(10, 0){{\tiny 10}}}
\put(160, 80){\makebox(10, 0){{\tiny 10}}}
\put(130, 10){\makebox(10, 0){\boxfill}}
\put(170, 10){\makebox(10, 0){\boxfill}}
\put(115, 15){\makebox(10, 0){\boxfill}}
\put(155, 15){\makebox(10, 0){\boxfill}}
\put(140, 20){\makebox(10, 0){\boxfill}}
\put(125, 25){\makebox(10, 0){\boxfill}}
\put(165, 25){\makebox(10, 0){\boxfill}}
\put(110, 30){\makebox(10, 0){\boxfill}}
\put(150, 30){\makebox(10, 0){\boxfill}}
\put(135, 35){\makebox(10, 0){\boxfill}}
\put(175, 35){\makebox(10, 0){\boxfill}}
\put(120, 40){\makebox(10, 0){\boxfill}}
\put(160, 40){\makebox(10, 0){\boxfill}}
\put(105, 45){\makebox(10, 0){\boxfill}}
\put(145, 45){\makebox(10, 0){\boxfill}}
\put(130, 50){\makebox(10, 0){\boxfill}}
\put(170, 50){\makebox(10, 0){\boxfill}}
\put(115, 55){\makebox(10, 0){\boxfill}}
\put(155, 55){\makebox(10, 0){\boxfill}}
\put(140, 60){\makebox(10, 0){\boxfill}}
\put(125, 65){\makebox(10, 0){\boxfill}}
\put(165, 65){\makebox(10, 0){\boxfill}}
\put(110, 70){\makebox(10, 0){\boxfill}}
\put(150, 70){\makebox(10, 0){\boxfill}}
\end{picture}
\caption{Expanding a $7\times 7$ square to a $9\times 9$ square. The dominoes given by the $7\times 7$ square are marked with ``\boxfill''.}
\label{77to99}
\end{figure}
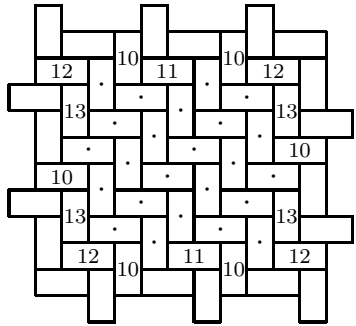

\medskip
By Lemma~\ref{tessalation}, the pairs of $\rho$ are either in the pattern $T_1$ or the pattern $T_2$, but none of those are pairings for Snaky. This concludes the proof of Theorem~\ref{nomatch}. \qed

\section{Torus games} \label{tori}

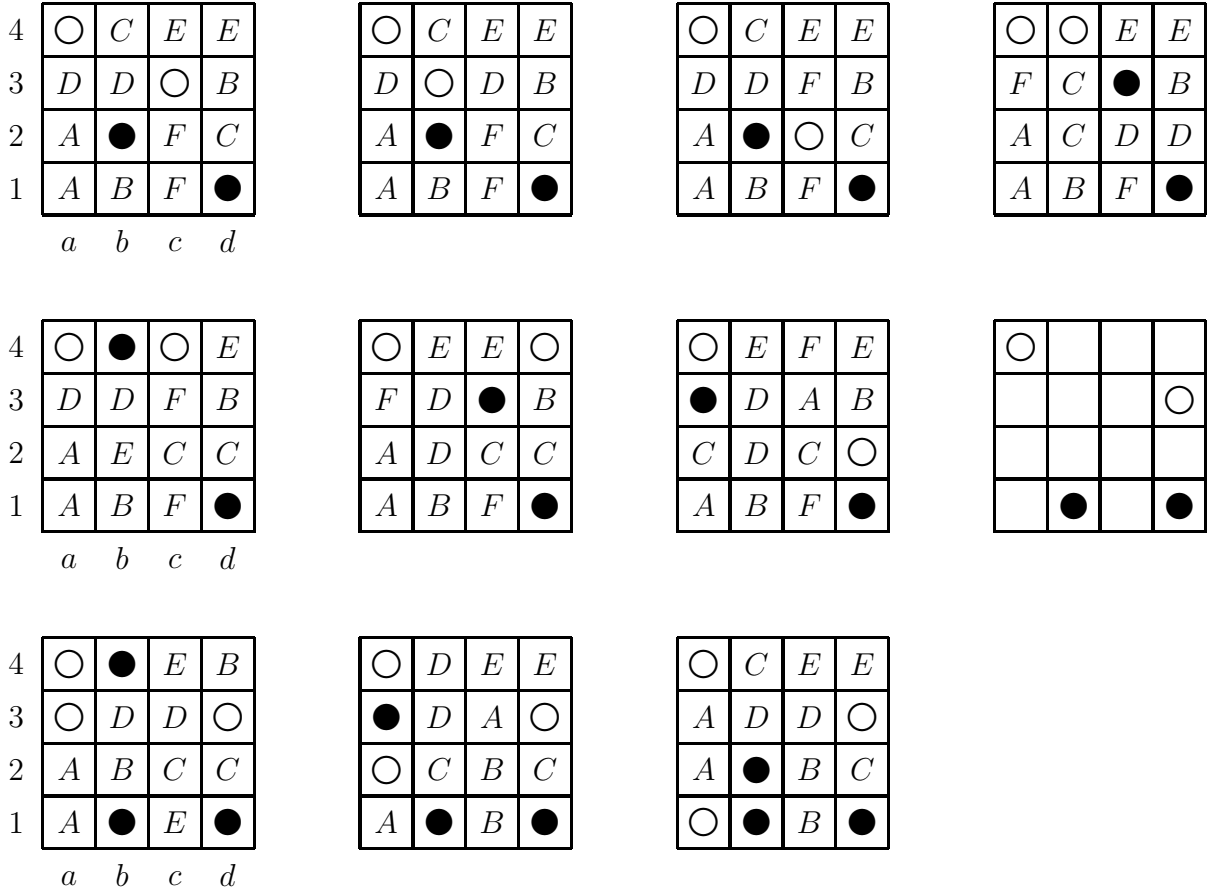
\begin{figure}[htbp]
\begin{picture}(400,400)(50,0)
\thicklines
%
\thicklines
\multiput(40,300)(0,20){5}{\put(0,0){\line(1,0){80}}}
\multiput(40,300)(20,0){5}{\put(0,0){\line(0,1){80}}}
\multiput(160,300)(0,20){5}{\put(0,0){\line(1,0){80}}}
\multiput(160,300)(20,0){5}{\put(0,0){\line(0,1){80}}}
\multiput(280,300)(0,20){5}{\put(0,0){\line(1,0){80}}}
\multiput(280,300)(20,0){5}{\put(0,0){\line(0,1){80}}}
\multiput(400,300)(0,20){5}{\put(0,0){\line(1,0){80}}}
\multiput(400,300)(20,0){5}{\put(0,0){\line(0,1){80}}}
%
\put(50, 370){\circle{10}}
\put(90, 350){\circle{10}}
\put(110, 310){\circle*{10}}
\put(70, 330){\circle*{10}}
\put(170, 370){\circle{10}}
\put(190, 350){\circle{10}}
\put(230, 310){\circle*{10}}
\put(190, 330){\circle*{10}}
\put(290, 370){\circle{10}}
\put(330, 330){\circle{10}}
\put(350, 310){\circle*{10}}
\put(310, 330){\circle*{10}}
\put(410, 370){\circle{10}}
\put(430, 370){\circle{10}}
\put(470, 310){\circle*{10}}
\put(450, 350){\circle*{10}}
%
\put(45, 289){\makebox(10, 0){$a$}}
\put(65, 290){\makebox(10, 0){$b$}}
\put(85, 289){\makebox(10, 0){$c$}}
\put(105,290){\makebox(10, 0){$d$}}
\put(25, 310){\makebox(10, 0){1}}
\put(25, 330){\makebox(10, 0){2}}
\put(25, 350){\makebox(10, 0){3}}
\put(25, 370){\makebox(10, 0){4}}
%
\put(45, 310){\makebox(10, 0){$A$}}
\put(45, 330){\makebox(10, 0){$A$}}
%
\put(65, 310){\makebox(10, 0){$B$}}
\put(105, 350){\makebox(10, 0){$B$}}
%
\put(65, 370){\makebox(10, 0){$C$}}
\put(105, 330){\makebox(10, 0){$C$}}
%
\put(45, 350){\makebox(10, 0){$D$}}
\put(65, 350){\makebox(10, 0){$D$}}
%
\put(85, 370){\makebox(10, 0){$E$}}
\put(105, 370){\makebox(10, 0){$E$}}
%
\put(85, 310){\makebox(10, 0){$F$}}
\put(85, 330){\makebox(10, 0){$F$}}
\put(165, 310){\makebox(10, 0){$A$}}
\put(165, 330){\makebox(10, 0){$A$}}
%
\put(185, 310){\makebox(10, 0){$B$}}
\put(225, 350){\makebox(10, 0){$B$}}
%
\put(185, 370){\makebox(10, 0){$C$}}
\put(225, 330){\makebox(10, 0){$C$}}
%
\put(165, 350){\makebox(10, 0){$D$}}
\put(205, 350){\makebox(10, 0){$D$}}
%
\put(205, 370){\makebox(10, 0){$E$}}
\put(225, 370){\makebox(10, 0){$E$}}
%
\put(205, 310){\makebox(10, 0){$F$}}
\put(205, 330){\makebox(10, 0){$F$}}
%
\put(285, 310){\makebox(10, 0){$A$}}
\put(285, 330){\makebox(10, 0){$A$}}
%
\put(305, 310){\makebox(10, 0){$B$}}
\put(345, 350){\makebox(10, 0){$B$}}
%
\put(305, 370){\makebox(10, 0){$C$}}
\put(345, 330){\makebox(10, 0){$C$}}
%
\put(285, 350){\makebox(10, 0){$D$}}
\put(305, 350){\makebox(10, 0){$D$}}
%
\put(325, 370){\makebox(10, 0){$E$}}
\put(345, 370){\makebox(10, 0){$E$}}
%
\put(325, 310){\makebox(10, 0){$F$}}
\put(325, 350){\makebox(10, 0){$F$}}
%
%
\put(405, 310){\makebox(10, 0){$A$}}
\put(405, 330){\makebox(10, 0){$A$}}
%
\put(425, 310){\makebox(10, 0){$B$}}
\put(465, 350){\makebox(10, 0){$B$}}
%
\put(425, 350){\makebox(10, 0){$C$}}
\put(425, 330){\makebox(10, 0){$C$}}
%
\put(445, 330){\makebox(10, 0){$D$}}
\put(465, 330){\makebox(10, 0){$D$}}
%
\put(445, 370){\makebox(10, 0){$E$}}
\put(465, 370){\makebox(10, 0){$E$}}
%
\put(445, 310){\makebox(10, 0){$F$}}
\put(405, 350){\makebox(10, 0){$F$}}
%
\thicklines
%
\thicklines
\multiput(40,180)(0,20){5}{\put(0,0){\line(1,0){80}}}
\multiput(40,180)(20,0){5}{\put(0,0){\line(0,1){80}}}
\multiput(160,180)(0,20){5}{\put(0,0){\line(1,0){80}}}
\multiput(160,180)(20,0){5}{\put(0,0){\line(0,1){80}}}
\multiput(280,180)(0,20){5}{\put(0,0){\line(1,0){80}}}
\multiput(280,180)(20,0){5}{\put(0,0){\line(0,1){80}}}
\multiput(400,180)(0,20){5}{\put(0,0){\line(1,0){80}}}
\multiput(400,180)(20,0){5}{\put(0,0){\line(0,1){80}}}
%
\put(50, 250){\circle{10}}
\put(90, 250){\circle{10}}
\put(110, 190){\circle*{10}}
\put(70, 250){\circle*{10}}
\put(170, 250){\circle{10}}
\put(230, 250){\circle{10}}
\put(230, 190){\circle*{10}}
\put(210, 230){\circle*{10}}
\put(290, 250){\circle{10}}
\put(350, 210){\circle{10}}
\put(350, 190){\circle*{10}}
\put(290, 230){\circle*{10}}
\put(410, 250){\circle{10}}
\put(470, 230){\circle{10}}
\put(470, 190){\circle*{10}}
\put(430, 190){\circle*{10}}
%
\put(45, 169){\makebox(10, 0){$a$}}
\put(65, 170){\makebox(10, 0){$b$}}
\put(85, 169){\makebox(10, 0){$c$}}
\put(105,170){\makebox(10, 0){$d$}}
\put(25, 190){\makebox(10, 0){1}}
\put(25, 210){\makebox(10, 0){2}}
\put(25, 230){\makebox(10, 0){3}}
\put(25, 250){\makebox(10, 0){4}}
%
\put(45, 190){\makebox(10, 0){$A$}}
\put(45, 210){\makebox(10, 0){$A$}}
%
\put(65, 190){\makebox(10, 0){$B$}}
\put(105, 230){\makebox(10, 0){$B$}}
%
\put(85, 210){\makebox(10, 0){$C$}}
\put(105, 210){\makebox(10, 0){$C$}}
%
\put(45, 230){\makebox(10, 0){$D$}}
\put(65, 230){\makebox(10, 0){$D$}}
%
\put(65, 210){\makebox(10, 0){$E$}}
\put(105, 250){\makebox(10, 0){$E$}}
%
\put(85, 190){\makebox(10, 0){$F$}}
\put(85, 230){\makebox(10, 0){$F$}}
%
\put(165, 190){\makebox(10, 0){$A$}}
\put(165, 210){\makebox(10, 0){$A$}}
%
\put(185, 190){\makebox(10, 0){$B$}}
\put(225, 230){\makebox(10, 0){$B$}}
%
\put(205, 210){\makebox(10, 0){$C$}}
\put(225, 210){\makebox(10, 0){$C$}}
%
\put(185, 210){\makebox(10, 0){$D$}}
\put(185, 230){\makebox(10, 0){$D$}}
%
\put(185, 250){\makebox(10, 0){$E$}}
\put(205, 250){\makebox(10, 0){$E$}}
%
\put(205, 190){\makebox(10, 0){$F$}}
\put(165, 230){\makebox(10, 0){$F$}}
%
\put(285, 190){\makebox(10, 0){$A$}}
\put(325, 230){\makebox(10, 0){$A$}}
%
\put(305, 190){\makebox(10, 0){$B$}}
\put(345, 230){\makebox(10, 0){$B$}}
%
\put(325, 210){\makebox(10, 0){$C$}}
\put(285, 210){\makebox(10, 0){$C$}}
%
\put(305, 210){\makebox(10, 0){$D$}}
\put(305, 230){\makebox(10, 0){$D$}}
%
\put(305, 250){\makebox(10, 0){$E$}}
\put(345, 250){\makebox(10, 0){$E$}}
%
\put(325, 190){\makebox(10, 0){$F$}}
\put(325, 250){\makebox(10, 0){$F$}}
%
\thicklines
\thicklines
\multiput(40,60)(0,20){5}{\put(0,0){\line(1,0){80}}}
\multiput(40,60)(20,0){5}{\put(0,0){\line(0,1){80}}}
\multiput(160,60)(0,20){5}{\put(0,0){\line(1,0){80}}}
\multiput(160,60)(20,0){5}{\put(0,0){\line(0,1){80}}}
\multiput(280,60)(0,20){5}{\put(0,0){\line(1,0){80}}}
\multiput(280,60)(20,0){5}{\put(0,0){\line(0,1){80}}}
%
\put(50, 130){\circle{10}}
\put(110, 110){\circle{10}}
\put(50, 110){\circle{10}}
\put(110, 70){\circle*{10}}
\put(70, 70){\circle*{10}}
\put(70, 130){\circle*{10}}
\put(170, 130){\circle{10}}
\put(230, 110){\circle{10}}
\put(170, 90){\circle{10}}
\put(230, 70){\circle*{10}}
\put(190, 70){\circle*{10}}
\put(170, 110){\circle*{10}}
\put(290, 130){\circle{10}}
\put(350, 110){\circle{10}}
\put(290, 70){\circle{10}}
\put(350, 70){\circle*{10}}
\put(310, 70){\circle*{10}}
\put(310, 90){\circle*{10}}
%
\put(45, 49){\makebox(10, 0){$a$}}
\put(65, 50){\makebox(10, 0){$b$}}
\put(85, 49){\makebox(10, 0){$c$}}
\put(105,50){\makebox(10, 0){$d$}}
\put(25, 70){\makebox(10, 0){1}}
\put(25, 90){\makebox(10, 0){2}}
\put(25, 110){\makebox(10, 0){3}}
\put(25, 130){\makebox(10, 0){4}}
%
\put(45, 70){\makebox(10, 0){$A$}}
\put(45, 90){\makebox(10, 0){$A$}}
%
\put(65, 90){\makebox(10, 0){$B$}}
\put(105, 130){\makebox(10, 0){$B$}}
%
\put(85, 90){\makebox(10, 0){$C$}}
\put(105, 90){\makebox(10, 0){$C$}}
%
\put(85, 110){\makebox(10, 0){$D$}}
\put(65, 110){\makebox(10, 0){$D$}}
%
\put(85, 130){\makebox(10, 0){$E$}}
\put(85, 70){\makebox(10, 0){$E$}}
%
\put(165, 70){\makebox(10, 0){$A$}}
\put(205, 110){\makebox(10, 0){$A$}}
%
\put(205, 90){\makebox(10, 0){$B$}}
\put(205, 70){\makebox(10, 0){$B$}}
%
\put(185, 90){\makebox(10, 0){$C$}}
\put(225, 90){\makebox(10, 0){$C$}}
%
\put(185, 130){\makebox(10, 0){$D$}}
\put(185, 110){\makebox(10, 0){$D$}}
%
\put(205, 130){\makebox(10, 0){$E$}}
\put(225, 130){\makebox(10, 0){$E$}}
\put(285, 90){\makebox(10, 0){$A$}}
\put(285, 110){\makebox(10, 0){$A$}}
%
\put(325, 90){\makebox(10, 0){$B$}}
\put(325, 70){\makebox(10, 0){$B$}}
%
\put(305, 130){\makebox(10, 0){$C$}}
\put(345, 90){\makebox(10, 0){$C$}}
%
\put(325, 110){\makebox(10, 0){$D$}}
\put(305, 110){\makebox(10, 0){$D$}}
%
\put(325, 130){\makebox(10, 0){$E$}}
\put(345, 130){\makebox(10, 0){$E$}}
%
\end{picture}
\caption{The pairings used by Picker in the game $4^2$.}
\label{torusz}
\end{figure}

To test Beck's paradigm from Conjecture~\ref{Bc}that Chooser-Picker and Picker-Chooser games are similar to Maker-Breaker games, we check the status of concrete games defined on the $4 \times 4$ torus. That is, we identify the opposite sides of the grid, and consider all lines of slopes $0$ and $\pm 1$ and size $4$ to be winning sets.  We denote the torus, along with those winning sets with the notation $4^2$. For the general definition of torus games, see \cite{beckcpc05}. We use a chess-like notation to refer to the elements of the board. We note that the hypergraph of winning sets on $4^2$ is not almost disjoint, see e.\ g.\ the two winning sets $\{a2,b1,c4,d3\}$ and $\{a4,b1,c2,d3\}$.  See Figure~\ref{torusz}. We consider four possible games on $4^2$: Maker-Maker, Maker-Breaker, Chooser-Picker and Picker-Chooser. According to \cite{beckcpc05}, the Maker-Maker version of $4^2$ is a draw, and, according to~\cite{CMP}, Picker wins the Chooser-Picker version. Here, we investigate the Maker-Breaker and the Picker-Chooser versions. In fact, the statement of the Maker-Breaker version implies the result for the Maker-Maker version, while the proof of it contains the proof of the Chooser-Picker version.

\begin{proposition} \label{MBtorus}
Breaker wins the Maker-Breaker version of the $4^2$ torus game.
\end{proposition}

\noindent {\bf Proof.} Using the symmetry of $4^2$, we may assume, without loss of generality, that Maker takes $a4$. Breaker's move will then be to take $d1$. Up to isomorphism, there are eight cases depending on the next move of Maker. The first element of the pair is Maker's move, while the second is Breaker's answer: 1. $(c3,b2)$, 2. $(b3,b2)$, 3. $(c2,b2)$, 4. $(b4,c3)$, 5. $(c4,b4)$, 6. $(d4,c3)$, 7. $(d2,a3)$ and 8. $(d3,b1)$.

In the first seven cases Breaker has winning pairing strategies. All eight cases are shown in the first two rows of Figure~\ref{torusz} and the pairs appear under the labels $A$, $B$, $C$, $D$, and $E$. We leave it to the reader to check that the pairs block all 16 winning sets.

In the eighth case Breaker does not have pairing strategy, but the game reduces to one of the seven prior cases unless Maker plays $a3$, $a2$ or $a1$ in the third step of the game. In that case, Breaker plays $b4$, $a3$ or $b2$, respectively, and wins by the pairing strategy shown in the third row of Figure~\ref{torusz}. \qed
\medskip

Note that in the Chooser-Picker version of the game $4^2$, Picker can achieve a position isomorphic to Case 1.  That is, Picker wins.

If Conjecture~\ref{Bc} were true, then Breaker has an easier job in the Maker-Breaker version than Chooser has in the Picker-Chooser game. For the $4 \times 4$ torus the outcome of these games are the same, although this is much harder to prove.

\begin{proposition}\label{PCtorus}
Chooser wins the Picker-Chooser version of $4^2$, the $4 \times 4$ torus game.
\end{proposition}

\noindent {\bf Sketch of the proof.} The full proof needs a lengthy exhaustive case analysis. However, some branches of the game tree may be cut by the following result of Beck~\cite{Beck5}: Chooser wins the Picker-Chooser game on ${\mathcal H}$ if $T({\mathcal H}):=\sum_{A \in E({\mathcal H})} 2^{-|A|} <1$.

In our case, $T({\mathcal H})=16\times 2^{-4}=1$, which just falls short. Instead we use a similar method using so-called {\em potential functions}.\footnote{Decided to name it here. --RM} We assign weights to each edge at the $i^{\:\rm th}$ stage such that $w_i(A)=0$ if Chooser has taken an element of $A$, otherwise it is $2^{-f(A)}$, where $f(A)$ is the number of untaken elements of $A$. The weight of a vertex $x$ is $w_i(x)=\sum_{x \in A} w_i(A)$, while the total weight is $w_i:=\sum_{A \in E({\mathcal H})} w_i(A)$.

Note that Picker wins if and only if both $w_8 \geq 1$ and $w_0=T({\mathcal H})=1$. When a pair $(x,y)$ is offered, Chooser can always take the one with larger weight, which results in a non-increasing total weight. In fact, if the weights of $x$ and $y$ differ or both $x$ and $y$ are elements of an $A$ of positive weight, then the total weight strictly decreases.

In order to have any possibility of winning, Picker has to select $x$ and $y$ of equal weights and no edge of positive weight containing both. By the symmetries of the board, we may assume Picker gets $a4$ and Chooser gets $c3$ in the first round. After that, Picker has only pairs $(x,y)$ that do not result in a loss for Picker: $(b4, d3)$, $(a3, c4)$, $(b3, d4)$, $(a3, b3)$, $(a3, d3)$, $(b3, d3)$, $(a1, b2)$ and $(a1, d2)$, see Figure~\ref{pctorusz}. The letter P [C] designates the vertex taken by Picker [Chooser] in the first step, the numbers are the weights of the vertices.

\begin{figure}[htbp]
\begin{picture}(400,100)(50,0)
\thicklines
\thicklines
%
\multiput(220,10)(0,20){5}{\put(0,0){\line(1,0){80}}}
\multiput(220,10)(20,0){5}{\put(0,0){\line(0,1){80}}}
%
%
%
\put(226,76){\mbox{P}}
\put(266,56){\mbox{C}}
\put(225, 0){\makebox(10, 0){$a$}}
\put(245, 2){\makebox(10, 0){$b$}}
\put(265, 0){\makebox(10, 0){$c$}}
\put(285,2){\makebox(10, 0){$d$}}
\put(205, 20){\makebox(10, 0){1}}
\put(205, 40){\makebox(10, 0){2}}
\put(205, 60){\makebox(10, 0){3}}
\put(205, 80){\makebox(10, 0){4}}
\put(225, 20){\makebox(10, 0){$\frac{3}{16}$}}
\put(245, 20){\makebox(10, 0){$\frac{5}{16}$}}
\put(265, 20){\makebox(10, 0){$\frac{3}{16}$}}
\put(285, 20){\makebox(10, 0){$\frac{5}{16}$}}
\put(225, 40){\makebox(10, 0){$\frac{5}{16}$}}
\put(245, 40){\makebox(10, 0){$\frac{3}{16}$}}
\put(265, 40){\makebox(10, 0){$\frac{5}{16}$}}
\put(285, 40){\makebox(10, 0){$\frac{3}{16}$}}
\put(225, 60){\makebox(10, 0){$\frac{4}{16}$}}
\put(245, 60){\makebox(10, 0){$\frac{4}{16}$}}
\put(285, 60){\makebox(10, 0){$\frac{4}{16}$}}
\put(245, 80){\makebox(10, 0){$\frac{4}{16}$}}
\put(265, 80){\makebox(10, 0){$\frac{4}{16}$}}
\put(285, 80){\makebox(10, 0){$\frac{4}{16}$}}
\end{picture}
\caption{The beginning of the Picker-Chooser $4^2$ game.}
\label{pctorusz}
\end{figure}

The rest of the proof is similar to that of the prior step: one needs to check that Chooser has winning strategy for each of the eight nontrivial responses of Picker. We omit the details. \qed

\section{Thanks}
We thank anonymous referee for uncovering two flaws in a previous version of the manuscript and also for the numerous suggestions which improved the presentation of the paper.

\end{document}